\documentclass[11pt]{article}

\usepackage{graphicx}
\usepackage{amsmath}
\usepackage{amsfonts}
\usepackage{amssymb}
\usepackage{amsthm}
\usepackage{eucal}
\usepackage{nccmath}
\usepackage{array}
\usepackage{bbm}
\usepackage[margin=.75in]{geometry}
\usepackage{mathrsfs}
\usepackage{stmaryrd}
\usepackage{xfrac}
\usepackage[affil-it]{authblk}
\usepackage[font={small,it}]{caption}

\newtheorem{theorem}{Theorem}

\newtheorem{conjecture}[theorem]{Conjecture}

\newcommand{\half}[0]{\medmath{\frac{1}{2}}}

\newcommand{\wt}[1]{\widetilde{#1}}

\newcommand{\E}[1]{\mathbb{E}\left[ #1 \right]}

\newcommand{\tr}[0]{\,\text{tr}\,}

\renewcommand{\d}[0]{\mathrm{d}}

\renewcommand{\half}[0]{\tfrac{1}{2}}

\begin{document}
\title{Fine structure of moments of the KMK transform\\ of the Poissonized Plancharel measure}
\author{Patrick Waters}
\affil{Temple University}
\maketitle
\abstract{We consider asymptotics behavior of Poissonized Plancharel measures as the poissonization parameter $n$ goes to infinity.  Recently Moll proved a convergent series expansion for statistics of a measure $\mu_\lambda$ which is the Kerov-Markov-Krein transform of the signed measure on corners a Jack-random partition $\lambda$.  The measure $\mu_\lambda$ is of interest because it behaves in some ways like the empirical measure on eigenvalues of a GUE-random matrix.  
We prove for the Poissonized Plancharel case that the large $n$ series for moments of $\mu_\lambda$ have a recursive structure as rational expressions in the generating function for Catalan numbers.  We discuss the analogy between our result and the fine structure of moments of the GUE.}

\section{Introduction}
Parallels between the theories of random matrices and random partitions are well known; for a general introduction see \cite{OkNotes}.
In particular, it is known that there is a strong analogy between the Gaussian Unitary Ensemble (GUE) and the Poissonized Plancharel (PP) measure:
\begin{alignat}{2}
dP_n(M)=& Z_{n}^{-1}\exp\left(-n \tr \left(\half M^2 \right)\right)\, \d M, & \qquad \qquad & \text{GUE}(n)    \label{GUE}  \\
P_n (\lambda) = & \frac{n^{|\lambda|}}{e^n|\lambda|! }   \times \frac{(\dim \lambda)^2}{|\lambda|!} . & \qquad \qquad & \text{PP}(n)  \label{PP}
\end{alignat}
Here $M$ is an $n\times n$ hermitian matrix, and $\d M$ is Lebesgue measure on all its real and imaginary degrees of freedom.  
In equation (\ref{PP}) one can see that the number of squares $|\lambda|$ of a PP$(n)$ random partition has a Poisson distribution with parameter $n$; and given $|\lambda|$, the partition $ \lambda$ is then drawn from the corresponding Plancharel measure.  The quantity $\dim \lambda$ is the dimension of the irreducible representation of the symmetric group $S(|\lambda|)$ corresponding to the partition $\lambda$.

An example of the analogy between the measures PP$(n)$ and GUE$(n)$ is that correlation functions for both measures exhibit sine kernel universality ``in the bulk", and Airy kernel universality ``at the edge".
Another parallel is that moments of the GUE can be computed by enumerating maps, which are a kind of graph embedded in a Riemann surface; wheras moments of Plancharel measures can be expressed in terms of numbers of branched coverings of the Riemann sphere satisfying some conditions.
Okounkov used this connection to show that (when correctly rescaled) the largest part of a PP$(n)$ random partition has asymptotically the same Tracy-Widom distribution as the largest eigenvalue of a GUE$(n)$ random matrix \cite{Okounkov}.

In this paper we study Kerov-Markov-Krein transform of the measure PP$(n)$, which we now define.
The profile of a partition is the upper boundary of its Young diagram when rotated as in figure 1.  
\begin{figure} \centering
\includegraphics[width=.5\textwidth]{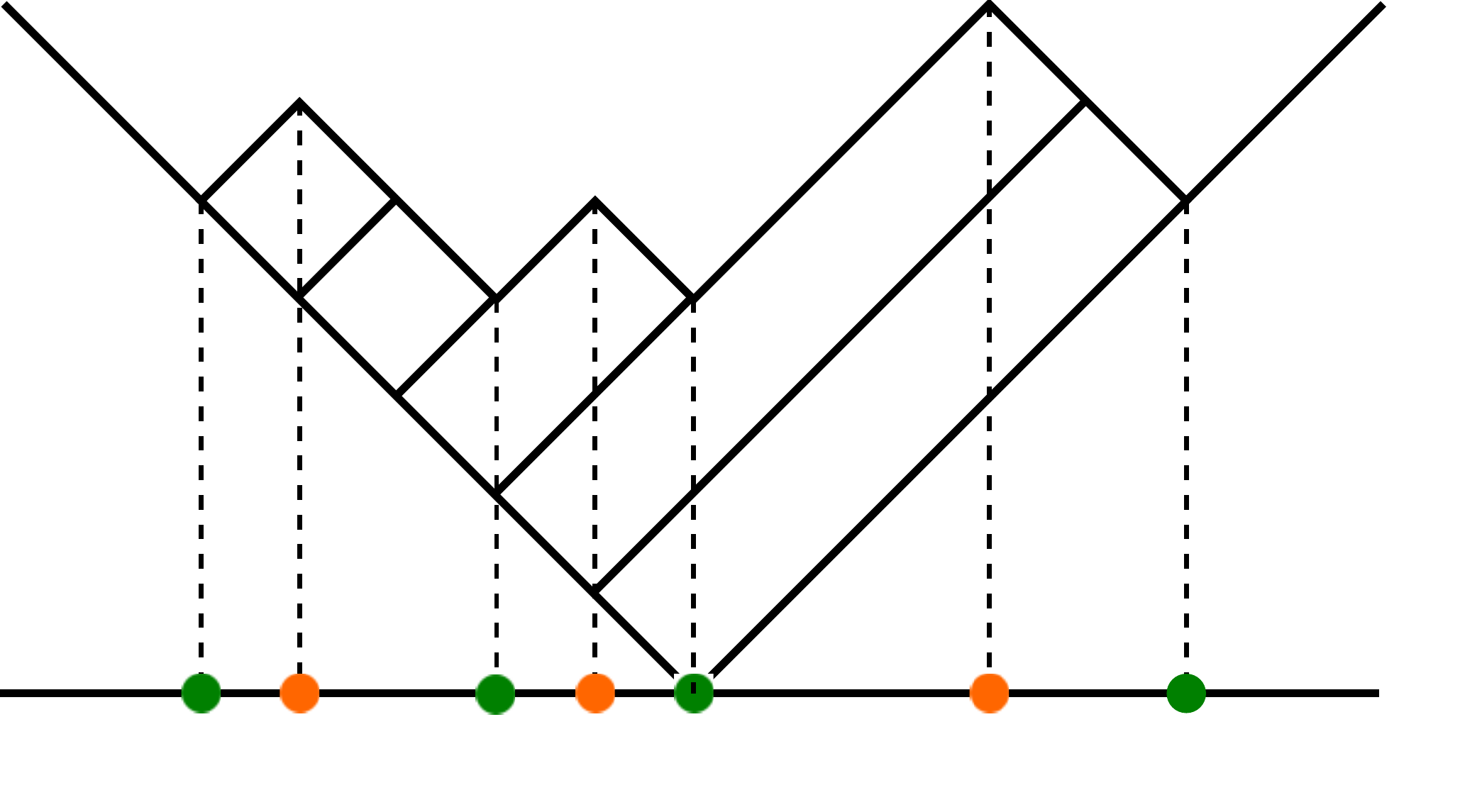}
\caption{Profile diagram for a partition $\lambda=5,5,2,1,1$.}  
\end{figure}
The diagram is drawn on a regular square lattice with side lengths $1/\sqrt{n}$, and thus in the scaling limit $n\rightarrow\infty$ has finite area.
Indeed, the Poissonized Plancharel measure has a well known limit shape \cite{LS77}.
The profile has upper and lower corners at horizontal coordinates $c^{\uparrow}_k(\lambda)$ and $c^{\downarrow}_k(\lambda)$ respectively.
It is easy to see that the upper corners of $\lambda$ are the elements of $\{(\lambda_i -i)/\sqrt{2n} \} \backslash \{ (\lambda_i-i-1)/\sqrt{2n}\}$, and a similar formula gives the locations of the lower corners.
The signed corner measure of the partition is the following superposition of Dirac masses:
\begin{align}
\sigma_\lambda =&\sum_k \delta_{c^{\uparrow}_k(\lambda)}- \delta_{c^{\downarrow}_k(\lambda)}.\nonumber
\end{align}
Following the recent paper of Moll \cite{MollQBO}, we define a transformed measure $\mu_\lambda$ by
\begin{align}
\exp\left(-\frac{1}{2} \int \log(x-s)\,\d \sigma_\lambda(s)\right)
=& \int \frac{d\mu_\lambda(s)}{x-s} \nonumber.
\end{align}
The new measure $\mu_\lambda$ is a probability measure.
In summary, the Poissonized Plancharel measure on partitions induces a distribution on probability measures $\mu_\lambda$ on $\mathbb{R}$, which is called the KMK transform of PP$(n)$.

The random measure $\mu_\lambda$ constructed above will be our principal object of study.  
We will exhibit a particular parallel between $\mu_\lambda$, and the empirical measure $\frac{1}{n} \sum_{i=1}^n \delta_{\eta_i}$ on the eigenvalues $\eta_i$ of a GUE$(n)$ random matrix.
Our starting point is the well known fact that the large $n$ limits of both of these random measure converge weakly to a point mass at the semicircle distribution
\begin{align}
\frac{1}{2\pi} \sqrt{4-s^2}\,  \mathbbm{1}_{[-2,2]}(s) \, \d s. \label{Wigner}
\end{align}
This motivates our interest in the measure $\mu_\lambda$ instead of (for example) an empirical measure on the downward steps of the partition profile or the signed corner measure.
Convergence to the semicircle measure (\ref{Wigner}) implies the following asymptotics of moment sequences:
\begin{align}
[n^0]\mathbb{E}_{PP(n)}\left[ \int \frac{\d \mu_\lambda(\eta)}{1- x\eta} \right] =&  c(x^2) = \frac{1-\sqrt{1-4x^2}}{2x^2} 
= [n^0] \int \frac{\rho^{(1)}_{\text{GUE}(n)}(\eta)}{1-x\eta}\d \lambda \label{momentLimit}.
\end{align}
Here the notation $[n^0]$ means the zeroth order term in a series expansion in powers of $n^{-1}$ for large $n$; this is of course just the limit as $n\rightarrow \infty$.  
It was proved in \cite{MollQBO} that the expression on the left hand side has a convergent series in powers of $n^{-1}$. 
The function $c(x^2)$ is a generating function for Catalan numbers
\begin{align}
c(x^2)= \sum_{k\geq 0} C_k x^{2k} = 1+x^2+2x^4+5x^6+14x^8+42x^{10}+\ldots,
\end{align}
and appears in the solutions of a variety of combinatorial problems.

On the right hand side of (\ref{momentLimit}), $\rho^{(1)}_{\text{GUE}(n)}$ is the one point correlation function for the GUE, defined as the averaged empirical measure on eigenvalues.
It is well known that the expression on the right hand side has a convergent series in $n$ if $x$ is nonreal,
in fact (as we will discuss in more detail in section \ref{GUEsection}) there exist integers $\ell_g(k)$ such that for $g\geq 1$
\begin{align}
[n^{-2g}] \int   \frac{\rho^{(1)}_{\text{GUE}(n)}(\lambda)\,}{1- x\lambda} \d \lambda =& \frac{c}{(2-c)^2} \sum_{k=2g}^{6g-3} \ell_g(k) \left(\frac{c-1}{2-c}\right)^k ,\qquad c=c(x^2).  \label{FineStructureGUE}
\end{align}
The notation $c=c(x^2)$ will be in force everywhere below.
This of course begs the question of whether some analogous result holds for the measures $\mu_\lambda$; our main theorem is the following affirmative answer:
\begin{theorem} \label{mainThm}
For each $g\in \mathbb{N}$ there rational numbers $\theta_{g}(k),\,k=0,\ldots ,2g$ such that
\begin{align}
[n^{-g}] \mathbb{E} \left[\int \frac{d\mu_\lambda(\eta)}{1-x\eta} \right] =&
\frac{c}{(2-c)^g} \sum_{k=g+1}^{3g-1} \theta_g(k) \left(\frac{c-1}{2-c}\right)^k  .\nonumber
\end{align}
\end{theorem}
We will use the notation $\Phi_g(c)$ for the functions on the right hand side of theorem \ref{mainThm}. 
Using a computed algebra system, we have computed the first few of the functions $\Phi_g$.  They are given by the following table of coefficients:
\begin{align}
\begin{array}{ll|llllllllll}
&& &&&&&k&&&& \\
&&2&3&4&5&6&7&8&9&10&11 \\ \hline
&1&1 & & & & & & & & &   \\
g&2& &1 & 14& 15 & & & & & &   \\
&3& & & 1 & 64 & 565 & 1122 & 630 & & &   \\
&4& & & &1 &222 & 5820&42500 & 110670 &118740 & 45045\\
\end{array} \nonumber
\end{align}
In section (\ref{GUEsection}) we will review the origins of formula (\ref{FineStructureGUE}).  This will lead us to make the following conjecture:
\begin{conjecture} \label{conj1}
The functions on the right hand side in theorem \ref{mainThm} are generating functions for some family of branched coverings of the sphere,
and the coefficients $\theta_g(k)$ enumerate a subfamily of them.
\end{conjecture}

\section{Fine structure of moments for GUE} \label{GUEsection}

In this section we give a quick review of the theory leading to equation (\ref{FineStructureGUE}) and conjecture \ref{conj1}.

It is well known that
\begin{align}
\int \lambda^{2k} \rho^{(1)}_{\text{GUE}(n)}(\lambda)\, \d \lambda =& \sum_{g\geq 0} \# \left\{ \substack{\text{Maps of genus $g$ with}\\ \text{one face and $k$ edges} }\right\} n^{-2g}. \label{oneFace1}
\end{align}
For each $k$ the sum on the right hand side is finite.
The classic paper on connection between map enumeration and matrix integrals was the work of Bessis, Itzykson and Zuber \cite{BIZ}, but (\ref{oneFace1}) was probably known earlier.
The right hand side was analyzed by Harer and Zagier in \cite{HZ86}.
Let $\epsilon_g(k)$ be the number of one faced maps of genus $g$ and $k$ edges, and $\ell_g(k)$ be the number of such maps with the further restriction that the map has no vertices of valences $1$ or $2$ (the valence of a vertex is the number of edges meeting there).
Define the generating functions
\begin{align}
E_g(x)=&\sum_{j} \epsilon_g(j) x^j, \qquad L_g(x) =\ell_g(j) x^j.
\end{align}
As above let $c=c(x^2)$.  It is shown in section $2$ of \cite{HZ86} that
\begin{align}
E_g  =& \frac{c}{(2-c)^2} L_g\left( \frac{c-1}{2-c} \right),\qquad g>0.  \label{elformula}
\end{align}
Formula (\ref{elformula}) is not valid for $g=0$, but in that case $E_0=c$.
For one faced maps with vertices of valence at least $3$ it follows from Euler's formula $V-E+F=2-2g$ that the number of edges is at most $6g-3$; and thus the generating functions $L_g$ are actually polynomials.
Furthermore by definition the faces of a map must be homeomorphic to discs; or in other words a handle may not be contained within a face.
Thus at least two edges are required for each handle, so it follows that $L_g$ is divisible by $x^{2g}$.
This establishes the formula (\ref{FineStructureGUE}), and also gives an interpretation of the coefficients $\ell_g(k)$ in that formula as map counts.

As we mentioned above, it was used to great advantage in \cite{Ok00} that moments of Plancharel measures can be expressed in terms of counts of branched coverings of the Riemann sphere.
Specifically, let P$(n)$ be the Plancharel measure on partitions of $n$, and $U(\lambda)$ be the set of indices $i$ such that $\lambda_i -i$ is a lower corner of $\lambda$, then
\begin{align}
\mathbb{E}_{\text{P}(n)}\left[\sum_{i\in U(\lambda)} \delta_i(\lambda) (\lambda_i - i)^k \right] =& \# \bigg\{ \text{solutions to } (1) =\prod_{j=1}^k (1\, m_j) \text{ with each }m_j=2\ldots n \bigg\}, \label{OkMoment}\\
\text{where}\qquad \delta_i(\lambda) =& \frac{\lambda_i + \text{len}(\lambda)-i}{|\lambda|} \prod_{j:\; j\neq i} \left( 1- \frac{1}{\lambda_i-\lambda_j+j-i}\right).\nonumber
\end{align}
It is shown in \cite{Ok00} that the factors $\delta_i(\lambda)$ are approximately $1/\sqrt{n}$ with high probability, and thus the left hand side is approximately a kind of moment.
The equation on the right hand side is over permutations, and the solutions are in one to one correspondence with equivalence classes brached coverings of the sphere with the following monodromy specifications.  Let $p_1,\ldots , p_k \in S^2$ be the branch points of the covering and $\gamma_1,\ldots ,\gamma_k$ be small loops around them from a common basepoint $q$ with preimages $Q_1,\ldots ,Q_n$; then the conditions are that the monodromy actions of the $\gamma_i$ are transpositions, each $\gamma_i$ transposing $Q_1$ with some other $Q_{m(i)}$.

Since, at least from the perspective of moments, the measure $\mu_\lambda$ seems to bear a closer relation to the GUE than the partition corner measure appearing in (\ref{OkMoment}), it may be hoped that there is an analogous moment formula without any factor like the $\delta_i(\lambda)$.  This is the basis for conjecture \ref{conj1}.

We now describe another behavior of moments of GUE which we find analogous behavior of the KMK transform of the Poissonized Plancharel measure.
Kerov discovered that GUE moments have a combinatorial formula in terms of rook placements on Young diagrams \cite{KerovRooks}, which we now explain.
A placement of $k$ rooks on a Young diagram is a subset of $k$ squares such that no two squares have the same row or column.
Kerov's formula\footnote{
Our notation differs from \cite{KerovRooks} because our eigenvalue axis has been rescaled by $\sqrt{n}$ and our one-point correlation function is normalized to be a probability measure.
} is
\begin{align}
\int x^{2k} \rho^{(1)}_n(x)\, \d x =& \frac{1}{n^{k+1}} \sum_{s=0}^{n-1}\# R(s,k),
\end{align}
where $R(s,k)$ is set of all $k$-rook placements on Young diagrams with $k$ parts and each part in $\{s,s+1,\ldots ,s+k\}$.
Given a rook placement in $Y(s,k)$, let $c$ be the number of rooks in the first $s$ rows.  
Deleting all of the first $s$ rows that do not contain a rook gives a map
$
\bigcup_{s=0}^{n-1} R_\ell(s,k) \rightarrow R'(\ell,k), 
$
where $R'(\ell,k)$ is the subset of $R(\ell,k)$ such that exactly $\ell$ rooks are placed in the first $\ell$ rows.
The number of preimages of a point under this mapping can be computed using the formula $\sum_{s=0}^{n-1} s^{\underline{\ell}} =n^{\underline{\ell+1}}/(\ell+1)$, where $s^{\underline{\ell}}=s(s-1)\ldots (s-\ell+1)$.
This results in the formula
\begin{align}
\int x^{2k} \rho^{(1)}_n(x)\, \d x =& \frac{1}{n^{k+1}} \sum_{\ell=0}^k \binom{n}{\ell+1} \# R'(\ell,k). \label{kerov1}
\end{align}
Thus moments of the GUE can be expressed using a series in $n$ with coefficients given in terms of numbers rook placements on Young diagrams.
In section \ref{partitionSection} we give an analogous rook configuration formula (\ref{rooks4}) for moments of the measures $\mu_\lambda$.

\section{Moll's expansion}
Recently Moll discovered a convergent expansion for ``transformed statistics" of Jack measures which can be seen as an analogue of the topological expansion for random matrices \cite{MollQBO}.
A Jack measure is determined by a potential function $V$ and parameters $\epsilon=1/\sqrt{n}$ and $\beta$, and the Poissonized Plancharel measure we consider here is the special case $V(w)=w+w^{-1} ,\, \beta=2$.  
The ``all orders expansion" in \cite{MollQBO} holds for multivariate statistics of the measure $\d\mu_\lambda$ induced by an arbitrary Jack measure on $\lambda$,
but for brevity we will only restate the special case of this expansion for linear statistics and the Poissonized Plancharel measure.
Define a differential operator $V_{-1}=n^{-1}\partial_{V_1}$.
Then $V_1$ and $V_{-1}$ satisfy the canonical commutation relation\footnote{
In the case of a general Jack measure, there are a sequence of parameters $V_k$ and operators $V_{-k}$ which should be thought of as an algebra of Bosonic creation and annihilation operators.  The Jack polynomials can by defined in a simple way in terms of vacuum expectations of vertex operators acting on states in a Fermionic Fock space.  We refer the interested reader to \cite{MJD} for the special case of Schur measures; it is easy to extend the construction given there to general $\beta$ (i.e. a Jack measure).
}
$[V_{-1},V_1]=n^{-1}$,
but they commute with $\wt{V}_1$.
Define a function $\Pi$ and a semi-infinite matrix $L$ by
\begin{align}
\Pi =&\exp(nV_1 \wt{V}_1),\qquad
L_{ij}=\begin{cases} V_{\pm 1} & \text{ if }j-i =\pm 1\\
0 & \text{ else.}\end{cases} \nonumber 
\end{align}
The indices of $L$ are taken from the set $\{0,1,2,3,\ldots\}$.
The connection between the matrix $L$ and the corners of $\lambda$ was discovered by Nazarov and Sklyanin in \cite{NSQBO}, and it is a critical ingredient in the proof of the all orders expansion.
In the special case we consider here, the all orders expansion is
\begin{align}
\E{\int x^k \d \mu_{\lambda}(x)} =& \phi\left\{ \Pi^{-1} (L^{k})_{00} \Pi \right\}. \label{MollExp}
\end{align}
Here $\phi$ is the algebra homomorphism given by $V_1,\wt{V}_1\mapsto 1$, and not defined for $V_{-1}$.  Thus all $V_{-1}$'s must be resolved before the evaluation.  
This can be understood from the following example:
\begin{align}
\E{\int x^4 \d \mu_{\lambda}(x)}
=&\phi\left\{ \Pi^{-1}  (V_{1}^2 V_{-1}^2 + V_1 V_{-1}V_1 V_{-1}) \Pi \right\}  \nonumber \\
=&\phi\left\{ \Pi^{-1}  (V_{1}^2 V_{-1}^2 + V_1(V_1V_{-1}+n^{-1})V_{-1})  \Pi \right\} \nonumber  \\
=&\phi\left\{ \Pi^{-1}  (2V_{1}^2 \wt{V}_{1}^2 +n^{-1}V_{1}\wt{V}_1 ) \Pi \right\} \nonumber \\
=& 2+n^{-1}. \label{ex111}  
\end{align}
The procedure for evaluating higher moments is to use the commutation relation to move all $V_{-1}$'s to the right, where they act on $\Pi$ which has the effect of turning them into $\wt{V}_1$'s.

As discussed in \cite{MollQBO}, the terms that result from evaluating moments in this way correspond to a certain kind of marked lattice path.
It is easy to see that the terms of $(L^{2k})_{00}$ correspond to Dyck paths of length $2k$ constrained to always have nonnegative height.
These paths can be drawn on the lattice $\mathbb{Z}\times \mathbb{Z}$ as in figure 2, with up steps corresponding to factors $V_{1}$ and the down steps corresponding to factors $V_{-1}$.
In the example (\ref{ex111}) a new term is created when we apply the commutation relation $V_{-1}V_1 =V_1 V_{-1}+n^{-1}$.  
Thus terms in the moment calculation correspond not just to paths $p$, but to tuples $(p,m)$ where the ``path marking" $m$ is a choice of some downward steps $D=\{d_1,\ldots, d_k\}$, as many upward steps $U=\{u_1,\ldots ,u_k\}$, and a bijection $f:D\rightarrow U$ such that each $f(d_i)$ comes after $d_i$.  The marked path $(p,m)$ contributes $n^{-k}$ to $\mathbb{E}[\int x^r \, \d \mu_\lambda(x)]$ where $r$ is the length of $p$.
\begin{figure}
\centering
\includegraphics[width=.4\textwidth]{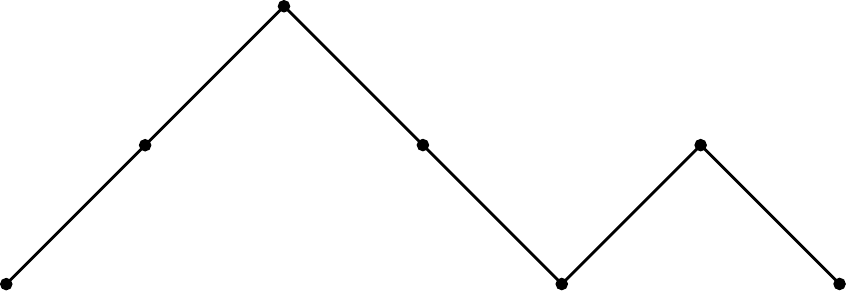}
\caption{This path contributes $V_{1}^2V_{-1}^2 V_1 V_{-1}$ to $(L^{6})_{00}$.}
\end{figure}

Let us consider the leading order terms in the expansion.  
This is the result one obtains by ignoring all the $n^{-1}$ commutator terms, and thus all nonnegative Dyck paths of length $2k$ make a contribution of $1$ to the $2k^{th}$ moment.
It is well known that Dyck paths are enumerated by Catalan numbers, that is
\begin{align}
[n^0] \E{\int x^{2k} \d \mu_{\lambda}(x)} 
=&\frac{1}{k+1}\binom{2k}{k}
 =[x^{2k}] \frac{1-\sqrt{1-4x^2}}{2x^2}.\nonumber 
\end{align}
Clearly all odd moments are zero because a lattice path must have an even number of steps if it is to begin and end at height zero.
Since Catalan numbers are the moments of the semicircle distribution (which has compact support), the above display shows that the limit shape of $\d \mu_\lambda$ is the semicircle as mentioned in the introduction.

\section{Lattice paths, partitions and rook placements}  \label{partitionSection}
We now address the problem of computing higher order terms of the expansion (\ref{MollExp}).
Our goal is to compute the functions
\begin{align}
\Phi_g(x)=& [n^{-g}] \mathbb{E}\left[\int \frac{\d\mu_\lambda(s)}{1-xs}\right]. \nonumber 
\end{align}
Let $P$ be the set of nonnegative Dyck paths starting at the point $(0,0)$, and ending at height zero.  Since the length of the paths is not constrained, $P$ is an infinite set.  The discussion in the previous section gives the following formula:  
\begin{align}
\Phi_g(x)
=& \sum_{p\in P} \# \left\{ \substack{ \text{markings of $p$ with} \\ \text{$g$ steps $\searrow$ and $g$ steps $\nearrow$} }\right\} x^{\text{length}(p)}.
\label{pathFormula1}
\end{align}
\begin{figure}
\centering
\includegraphics[width=.5\textwidth]{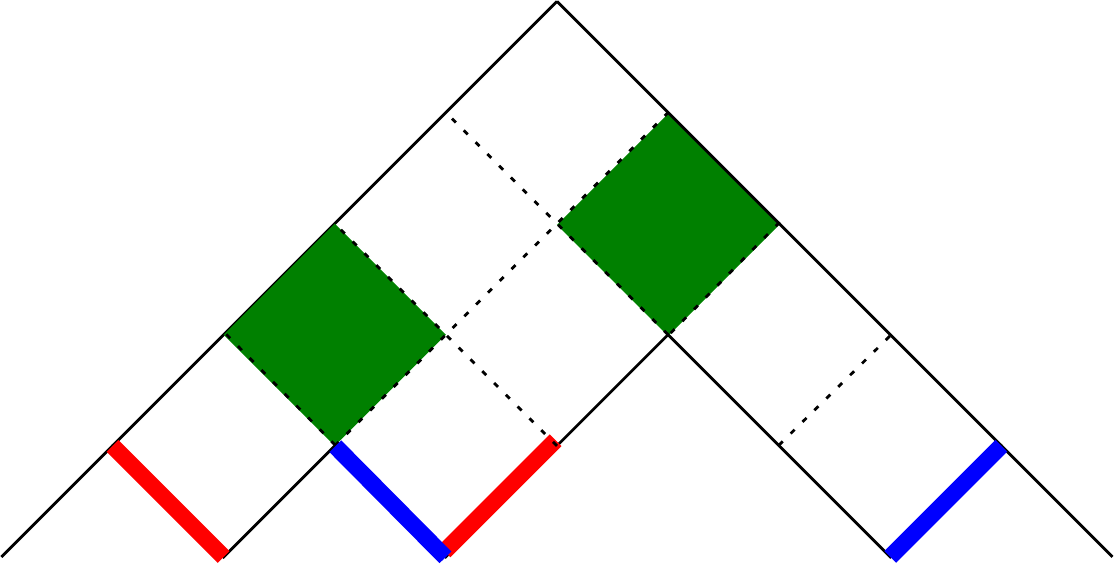}
\caption{A marked path and the corresponding rook placement}
\end{figure}
We now show that formula (\ref{pathFormula1}) can be interpreted in terms of rook placements on partitions.
Each path $p\in P$ of length $2k$ can be embedded in the triangle defined by points $(0,0),(k,k)$ and $(2k,0)$ because $p$ is nonnegative.  
The region of the triangle above $p$ is the Young diagram of a partition $\pi$, with parts corresponding to the downward steps of $p$.
Let $d_i(p)$ be the horizontal position just before the $i^{\text{th}}$ downward step of $p$, then the parts of $\pi$ are given by $\pi_i= k-1-p(i)+i$.

A marking of $p$ is indicated in figure 3 by giving a common color to steps that are paired.
Clearly a choice of one up step and one down step uniquely determines a box from our lattice, and it is easy to see that the box will be in the region above $p$ exactly if the downward step happens first.
Thus each marking of $p$ yields a subset $s$ of boxes of the Young diagram for $\pi$.
By construction, no two of these boxes can share a row or column, and thus $s$ is a rook placement on $\pi$.
Let $\Lambda_k$ be the set of partitions $\pi$ with parts satisfying the constraint $\pi_i \leq k-i$,
and $RC_r(k)$ be the set of $g$-rook configurations on partitions $\pi\in \Lambda_k$.
We have shown that
\begin{align}
  \mathbb{E}\left[\int x^{2k}\,\d\mu_\lambda(x)\right]
=& \sum_{r\geq 0}  \# RC_g(k)   n^{-g}. \label{rooks4}
\end{align}
This is a random partition analogue of Kerov's formula (\ref{kerov1}) for the GUE.

\section{Generating functions for paths}
To compute the correction terms using the formulas in the previous section, we need some formulas for Lattice path generating functions.
The formulas we present here must surely be well known, but we are unable to give a reference.
Let $P(i,j)$ be the set of nonnegative Dyck paths from $(0,0)$ to $(i,j)$.
Define
\begin{align}
F(x,y)=& \sum_{i,j\geq 0} \sum_{p\in P(i,j)} x^i y^j.\nonumber 
\end{align}
The function $F(x,y)$ can be expressed in terms of the generating function for Catalan numbers by decomposing each path its set of final returns to each height $(0,1,2,\ldots)$.
For any given path, there are a fixed finite number of heights to return to, but we simultaneously consider paths of all finite lengths.
It is thus easy to see that
\begin{align}
F(x,y)=\sum_{j\geq 0} (xy)^{j}c(x^2)^{j+1} = \frac{c(x^2)}{1-xyc(x^2)}.\nonumber 
\end{align}
We will also need the following slightly more general generating function.
Let $P(i,j_i,j_2)$ be the set of nonnegative Dyck paths from $(0,j_1)$ to $(i,j_2)$, then define
\begin{align}
G(x,y,z)=& \sum_{i,j_1,j_2\geq 0} \sum_{p\in P(i,j_1,j_2)} x^i y^{j_1} z^{j_2}.\nonumber 
\end{align}
By decomposing each path $p$ based on the last return of $p$ to its lowest height, one finds that
\begin{align}
G(x,y,z) =& \sum_{b\geq } y_{1}^b y_{2}^b F(x,y)F(x,z)/c(x^2) =\frac{F(x,y)F(x,z)}{c(x^2)(1-yz)}.\nonumber 
\end{align}

\section{Calculation of first correction term}
\begin{figure}
\centering
\includegraphics[width=.5\textwidth]{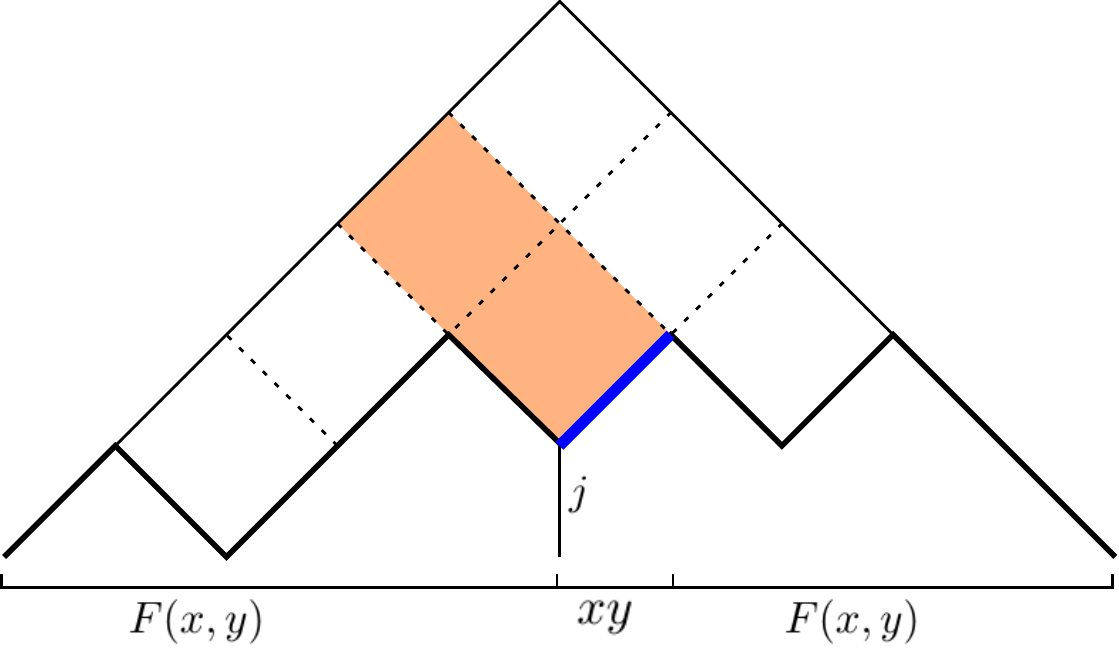}
\caption{Illustration of equation (\ref{asdf1})}
\end{figure}
Using the generating functions discussed in the previous section, we will now show that our rook placement formula (\ref{rooks4}) is equivalent to
\begin{align}
\Phi_1(x) =& \sum_{j\geq 0}  [y^{j+1}] xF(x,y) \times [y^j] \frac{1}{2} \left( x\partial_x -y \partial_y \right) F(x,y). \label{asdf1}
\end{align}
This formula can be understood as follows: we are to enumerate all ways of choosing a nonnegative Dyck path $p$ of length $2k$, and one box $t$ from the partition above $p$.
The tile $t$ will be chosen from a column $c$ determined by one of the up steps $s$ of the path.
There is a factor of $F(x,y)$ for the part of $p$ before step $s$.
The length of column $c$ is $(i-j)/2$ where $(i,j)$ is the point where step $s$ begins.
The differential operator in equation (\ref{asdf1}) acts on $F$ by multiplying the contribution of each path to $(i,j)$ by the desired factor $(i-j)/2$.
After step $s$, the path is once again arbitrary except that it must eventually terminate at height zero;
this part of $p$ is enumerated by another factor of $F(x,y)$. 
However there is the constraint that the path we take according to the second factor $F(x,y)$ must descend $j+1$ units, thus the series coefficients and sum over $j$.

The Catalan generating function satisfies:
\begin{align}
c(x)=1+x c(x)^2,\qquad c'(x)=\frac{c(x)^3}{ 2-c(x ) } .\nonumber 
\end{align}
Unless we specifically write $c(x)$, the notation $c=c(x^2)$ should be understood.
Using these identities, a straitforward calculation shows that
\begin{align}
\half (x\partial_x - y\partial_y )  F(x,y)=&\frac{ c(c-1)}{2-c}\frac{1}{(1- xyc)^2}.\nonumber 
\end{align}
Therefore we calculate
\begin{align}
\Phi_1(x)
&= \sum_{j\geq 0} [y^{j+1}]x\frac{c}{1-xyc} \times  [y^j] \frac{c(c-1)}{2-c} \frac{1}{(1-xy c)^2}  \nonumber \\
&=\sum_{j\geq 0} (xc)^{j+2} \times \frac{c(c-1)}{2-c} (j+1)(xc)^{j}\nonumber \\
&= \frac{c(c-1)^2}{(2-c)^3} .\label{phi1calc}
\end{align}
\section{Integral formula for $\Phi_g(x)$}

Before moving on to the general case, let us briefly consider $\Phi_2(x)$.  All of the considerations necessary for the general case will be observed here, and it is easier to explain the ideas in a concrete setting.
We will argue that
\begin{align}
\Phi_2(x)=&  \sum_{j_2\geq 0} [y^{j_2+1}] x F(x,y)  
 \times [y^{j_2}]  \left(\frac{1}{2}  x\partial_x -\frac{1}{2} y\partial_{y } -1\right)  H(x,y),\label{phi21} \\ 
H(x,y)=&\sum_{j_1\geq 0}  [z^{j_1+1}] x G(x,y ,z) 
\times [z^{j_1}]   \frac{1}{2}\left( x\partial_x -z \partial_{z}\right) F(x,z)  . \nonumber
\end{align}
By formula (\ref{rooks4}), $\Phi_2$ enumerates all configurations of two rooks on tiles $t_1,t_2$ of the partition above a Dyck path $p$, with each path weighted by a factor of $x^{\text{len}(p)}$.
We begin by selecting two up steps $s_1<s_2$ of $p$ (with $<$ meaning ``to the left of"); the steps $s_1,s_2$ determine two columns $c_1,c_2$ from which to choose $t_1,t_2$.
To the part of $p$ before $s_1$ there corresponds a factor of $F(x,z)$, and we operate on it with $\half(x\partial_x -z\partial_z)$ so that the contribution of each term will be multiplied be the number of tiles in $c_1$.
The factor of $G(x,y,z)$ enumerates the possibilities for the portion of $p$ between $s_1$ and $s_2$.
We once again apply an operator $\half(x\partial_x -y\partial_y)-1$, with the $-1$ to account for the restriction that $t_1,t_2$ cannot be taken from the same row.
The final factor of $F(x,y)$ describes the part of $p$ after $s_2$.

From the discussion above, the behavior of the general case is clear.
Define a sequence of operators $\mathcal{G}_k$ acting on functions $H(x,y)$ by
\begin{align}
(\mathcal{G}_k  H)(x,y)=& \sum_{j\geq 0} [z^{j+1}] x G(x,y,z) \times [z^j] \left(\frac{1}{2}x\partial_x -\frac{1}{2}z \partial_{z}  -k \right) H(x,z).\nonumber 
\end{align}
Then
\begin{align}
\Phi_g (x)=  [y^0] \mathcal{G}_{g-1} \mathcal{G}_{r-2}\ldots \mathcal{G}_0 \cdot F(x,y). \label{phiRFormula}
\end{align}

\section{Rational expressions for $\Phi_g$}
In this section we prove theorem \ref{mainThm}.
By a straitforward calculation we have
\begin{align}
[z^{j+1}]x G(x,y,z)= \frac{c-1}{1-x c y} \sum_{\ell=0}^{j+1}\left( \frac{y}{x c}\right)^\ell (x c)^{j}.
\end{align}
Letting $\mathcal{E}_r$ be the Euler operator $\frac{1}{2}(x\partial_x -z\partial_z)-r$, it follows that
\begin{align}
(\mathcal{G}_r H)(x,y) =& \sum_{j\geq 0} \frac{c-1}{1-x c y} \sum_{\ell=0}^{j+1}\left( \frac{y}{x c}\right)^\ell (x c)^{j} [z^j] \mathcal{E}_r H(x,z) \nonumber \\
=& \frac{c-1}{1-x c y}  \sum_{\ell\geq 0} \left( \frac{y}{x c}\right)^\ell  \sum_{j\geq \ell-1}(x c)^{j} [z^j] \mathcal{E}_r H(x,z). \label{d1}
\end{align}
Making a change of variables from $x$ to $c$, the Euler operator transforms by
\begin{align}
\frac{1}{2}(x \partial_x -y\partial_y) =\frac{c(c-1)}{2-c}\partial_c -\frac{1}{2}y\partial_y.\nonumber
\end{align}
It is easy to see that for any polynomial $P(c)$ and positive integers $s,k,m$, there are polynomials $\wt{P}_1$ and $\wt{P}_2$ such that
\begin{align}
\mathcal{E}_k  \frac{c(c-1)^s P(c)}{(2-c)^k(1-x c y)^{m}} =& \frac{c(c-1)^s \wt{P}_1(c)}{(2-c)^{k+2}(1-x c y)^{m}}+\frac{c(c-1)^{s} \wt{P}_2(c)}{(2-c)^{k+1}(1-x c y)^{m+1}}. \label{d2}
\end{align}
We now compute the parts of equation (\ref{d1}).  The functions of $c$ appearing in (\ref{d2}) have little effect on the rest of the calculation; what we are really interested is:
\begin{align}
  \sum_{j\geq \ell-1} (x c)^j [z^j] \frac{1}{(1-x c z)^{m+1}} 
=&  \sum_{j\geq \ell-1} \binom{j+m}{m}(c-1)^{2j}  
= \frac{1}{m!} \partial_{c}^{m} \sum_{j\geq \ell-1} (c-1)^{j+m} \nonumber \\
=&   \frac{1}{m!} \partial_{c}^m \frac{(c-1)^{\ell+m-1}}{2-c}  
=   (c-1)^\ell \sum_{i=0}^{m}  \frac{\wt{P}_{m,i}(c) \ell^{\underline{i}} }{ (2-c)^{1+m-i}}, \nonumber
\end{align}
for some coefficients $\wt{P}_{m,i}$ that are polynomials in $c$.
Here $\ell^{\underline{i}} =\ell(\ell-1)\ldots (\ell-i+1)$, and we chose this basis for polynomials in $\ell$ so that we will be able to recognize some multiple derivatives in the next step:
\begin{align}
 \sum_{\ell\geq 0} \left( \frac{y}{x c}\right)^\ell  \sum_{j\geq \ell-1}(x c)^{j} [z^j] \frac{1}{(1-x c y)^{m+1}}
=&\sum_{\ell\geq 0} \left( x c y\right)^\ell  \sum_{i=0}^{m}  \frac{\wt{P}_{m,i}(c) \ell^{\underline{i}} }{ (2-c)^{1+m-i}} \nonumber \\
=& \sum_{i=0}^m \frac{\wt{P}_{m,i}(c) }{ (2-c)^{1+m-i}} (x c y)^i \partial_{x c y}^i \sum_{\ell \geq 0} (x c y)^\ell \nonumber \\
=&\sum_{i=0}^m \frac{\wt{P}_{m,i}(c) }{ (2-c)^{1+m-i}} \frac{i!(x c y)^i}{(1-x c y)^{i+1}}. \label{d3}
\end{align}
Using formulas (\ref{d1},\ref{d2},\ref{d3}), it follows by induction that
\begin{align}
 \mathcal{G}_{r-1} \mathcal{G}_{r-2}\ldots \mathcal{G}_0 \cdot F(x,y) = & \sum_{i=0}^{2r-1} \frac{c(c-1)^{r} P_{r,i}(c)}{(2-c)^{4r-1-i}(1-x c y)^{2+i}}
\end{align}
for some polynomials $P_{r,i}(c)$ of degree $2r-1-i$.
Taking a series coefficient $[y^0]$ on both sides proves theorem \ref{mainThm}.

Using the formulas above and the help of a computer algebra system, we were able to calculate the table of coefficients $\theta_g(k)$ given in the introduction.
Given the number of errors in the first version of this paper, the reader would be right to question whether this table can be trusted.
However, we are confident that our calculations are correct because series expansions in $x$ of the formulas in our table match against numbers of rook placements that we computed separately using a different computer program.


\begin{thebibliography}{99}
\bibitem{BIZ}D. Bessis, C. Itzykson, J.B. Zuber, \emph{Quantum field theory techniques in graphical enumeration}, Adv. Appl. Math., 1, 109-157, 1980.
\bibitem{MJD} E. Date, M. Jimbo, T. Miwa, \emph{Solitons: differential equations, symmetries and infinite dimensional algebras}, Cambridge Tracts in Mathematics, 1999.
\bibitem{KerovRooks} S.V. Kerov, \emph{Rooks on Ferrers boards and matrix integrals}, J. Math. Sci., 96(5), 1999.
\bibitem{HZ86} J. Harer, D. Zagier, \emph{The Euler Characteristic of the moduli space of curves}, Invent. Math., 85 (1986) 457-485.
\bibitem{LS77}B.F. Logan and L.A. Shepp, \emph{A variational problem for random Young tableaux}, Adv. Math., 26, 206-222, 1977.
\bibitem{NSQBO}M.L. Nazarov, E.K. Sklyanin, \emph{Integrable hierarchy of the quantum Benjamin-Ono equation}, SIGMA, 9:78-92, 2013.
\bibitem{MollQBO}A. Moll, \emph{Random partitions and the Quantum Benjamin-Ono Hierarchy}, arXiv: 1508.03063, 2015.
\bibitem{Ok00}A. Okounkov, \emph{Random matrices and random permutations}, IMRN, 20, 1043-1095, 2000.
\bibitem{OkNotes}A. Okounkov, \emph{Symmetric functions and random partitions}, NATO Science Series, 74, 223-252, 2001.
\end{thebibliography}
\end{document}